\documentclass[11pt]{article}
\usepackage[utf8]{inputenc}
\usepackage{amsfonts, amsmath}
\usepackage{amssymb}
\usepackage{amsthm,amsmath,amscd}
\usepackage{enumerate}
\usepackage{url}
\usepackage[pdftex]{graphicx}
\usepackage[margin=25pt,font=small,labelfont=bf]{caption}
\usepackage{subfig}
\usepackage[numbers,square,sort&compress]{natbib}
\usepackage[colorlinks=true]{hyperref}
\newcommand{\defn}[1]{\textcolor{blue}{\emph{#1}}}
\newcommand*{\doi}[1]{doi: \href{https://dx.doi.org/#1}{\urlstyle{rm}\nolinkurl{#1}}}
\newcommand*{\arxiv}[1]{arXiv:  \href{https://arxiv.org/abs/#1}{\urlstyle{rm}\nolinkurl{#1}}}
\let\oldproofname=\proofname
\renewcommand{\proofname}{\rm\bf{\oldproofname}}

\newcommand{\RR}{\mathbb R}

\newcommand{\bna}{\begin{eqnarray}}
\newcommand{\ena}{\end{eqnarray}}
\newcommand{\ba}{\begin{eqnarray*}}
\newcommand{\ea}{\end{eqnarray*}}
\newcommand{\bs}[1]{}

\newcommand{\Rank}{{\rm Rank}}
\def\p{{\bf p}}

\def\0{{\bf 0}}

\let\oldv=\v
\def\v{{\bf v}}

\def\hp{{\hat{\p}}}
\def\hR{{\hat{R}}}

\usepackage[margin=1in]{geometry}
\begin{document}
\title{
Transverse Rigidity is Prestress Stability}

\author{
Steven J. Gortler
\and
Miranda Holmes-Cerfon
\and 
Louis Theran}
\date{}
\maketitle 

\begin{abstract}  
Recently, V. Alexandrov proposed an intriguing sufficient condition for rigidity,
which we will call ``transverse rigidity''. We show that
transverse rigidity
is actually equivalent to the known sufficient condition for rigidity called 
``prestress stability''. Indeed this leads to a novel
interpretation of the prestress condition.
\end{abstract}
In~\cite{va}, an  intriguing sufficient condition for rigidity was proposed,
which we will call ``transverse rigidity''. In this note, we show that, when it applies, 
transverse rigidity
is actually equivalent to the known sufficient condition for rigidity called 
``prestress stability'' (see ~\cite{Connelly-Whiteley}).
This equivalence was a bit of surprise to us, as these two 
conditions are based on different geometric intuitions, and some 
massaging is needed to see their equivalence.
Indeed the equivalence leads to a novel 
interpretation of the prestress condition and of second order rigidity.

Let $d \ge 1$ be a fixed dimension, and let $D:=\binom{d+1}{2}$.
Let $G$ be a graph that is generically rigid in $\RR^d$, with $n \ge d+1$ vertices and 
$m:=dn-D$ edges (it is isostatic).
Let $\hp$ denote a configuration of $n$ points in $\RR^d$ and
$(G,\hp)$ denote the associated bar-and-joint framework.
We assume that $\hp$ has a full $d$-dimensional span.

We denote by $\hR(\hp)$ the \defn{rigidity matrix} associated with $(G,\hp)$.
Its rows are indexed by edges of $G$: unordered distinct pairs $\{ij\}$. 
Its columns are indexed
by vertex $i \in \{1...n\}$ and spatial dimension $k \in \{1...d\}$. This matrix is $m$-by-$nd$.
Its non-zero entries are of the forms 
\ba
\hR(\hp)^{\{ij\}}_{ik} &=& (\hat{p}_i)_k - (\hat{p}_j)_k 
\ea

If $\hR(\hp)$ has rank $nd-D$, then we say that $(G,\hp)$ is \defn{infinitesimally rigid},
and this implies that $(G,\hp)$ is rigid. 
From our assumptions, all generic $(G,\hp)$ will be infinitesimally rigid. 
But what can we say for special $\hp$ where the rank of $\hR(\hp)$ drops?
We are interested in conditions that can be used to certify that 
some specific $(G,\hp)$ is still rigid.

To simplify the discussion going forward we
mod out the trivial motions around some specific $\hp$ by 
pinning down the framework.  It is well-known we can do this 
by fixing an appropriately chosen set of $D$  of the coordinates
in $\hp$ (see \cite{WW-pinning}).  The effect on the rigidity matrix is to remove the 
corresponding $D$ columns.  This leaves us with a pinned
``square rigidity matrix'' of size $(nd-D)$-by-$(nd-D)$, which 
we denote by (the unhatted) $R(\hp)$.

Going forward for now, let us suppose, that we have $\Rank(R(\hp))=nd-D-1$.
In this case there will be, up to scale,  a single row vector $\omega$ in the left kernel
of $R(\hp)$, this is called a \defn{equilibrium stress}.
(Note that this $\omega$ will also be the left kernel of $\hR(\hp)$.)
There will be, up to scale, a single column vector $\p'$ (of size $nd-D$) in its
right kernel, this is called an \defn{infinitesimal flex}.
We denote by $\hp'$, the vector (of size $nd$)
obtained by from $\p'$ by 
zero padding the pinned degrees of freedom.
The vector $\hp'$ is a non-trivial flex in the kernel
of $\hR(\hp)$.

With this notation, one way to certify rigidity is by testing for 
\defn{prestress stability}.
We check if
\bna
\label{eq:pss1}
\sum_{\{ij\}\in G} \omega_{\{ij\}} ||\hp'_i -\hp'_j||^2  \neq 0
\ena
If the condition \eqref{eq:pss1} holds, then $(G,\hp)$ is 
prestress stable and hence rigid
~\cite[Eq. (1), Props. 3.4.2 and 3.3.2]{Connelly-Whiteley}.
We note that there is another related notion called ``second-order'' rigidity, 
but it coincides with prestress stability when a framework 
has either one infinitesimal flex or one stress~\cite[Cor. 5.3.1]{Connelly-Whiteley}.

Alexandrov describes another rigidity test which we explain next.
The basic idea is that if there were a finite 
(pinned) flex $\hp(t)$ of $G$, starting at
$\hp(0)=\hp$, then during the flex the rigidity matrix would have to 
stay singular. Note that the singularity condition is
\ba
\det(R(\hp(t))) = 0
\ea
Differentiating wrt $t$ at $t=0$,  applying the chain rule, 
and then using the fact that the time derivative of $\hp(t)$ corresponds to  an infinitesimal flex,
we would get
\ba
d[\det(R(\hp))] \p' = 0
\ea
Here $d[f]$ is the $nd-D$ row vector representing the partials of $f$ with respect
to the $nd-D$ free variables in $\hp$.
Thus a sufficient condition for rigidity 
is
\bna
\label{eq:trans}
d[\det(R(\hp))] \p' \neq 0
\ena
Let us call this new condition \defn{transverse rigidity}, since it suggests the
idea that the potential direction for flexing, $\hp'$, lies transverse to the locus
of configurations with singular rigidity matrices.

\paragraph{Equivalence}
Here we show the equivalence of Equations (\ref{eq:pss1}) and (\ref{eq:trans}).
To begin with, 
Equation (\ref{eq:pss1}) can be rewritten 
~\cite[Eq. (1)]{Connelly-Whiteley}
as
\bna
\label{eq:fred}
\omega R(\hat{\p}') \p' \neq 0
\ena
We used the fact that $\hp'$ has zeros in the pinned coordinates to obtain an expression with $R$ instead of $\hR$. 
Note that in the term $R(\hp')$ we actually have a rigidity matrix built off
of the ``configuration'' $\hp'$.

Now our job is to show the equivalence of Equations (\ref{eq:trans}) and (\ref{eq:fred}).
We will do this by establishing  that
\bna
\omega R(\hat{\p}') \propto
d[\det(R(\hat{\p}))] 
\ena
(equality up to scale).

Let us start with the left side. From the definition of a rigidity matrix, 
we compute
\ba
[\omega R(\hat{\p}')]_{ik} 
&=& \sum_{j\in N(i)}  \omega_{\{ij\}}((\hat{p}'_i)_k - (\hat{p}'_j)_k)
\\
 &=& \sum_{j\in N(i)}  \omega_{\{ij\}}(p'_i)_k - \sum_{j\in N_k(i)}\omega_{\{ij\}}(p'_j)_k
\ea
Here $ik$ corresponds to an unpinned degree of freedom. 
$N(i)$ is the set of vertices that are neighbors of
$i$ in $G$. $N_k(i)$ is the set of neighbors $j$ of $i$
where the $k$th coordinate is not pinned

Now on to the right hand side.
We will compute $d[\det(R(\hp))]$ using the chain rule.

Given a square matrix $M$, we have
\ba
\frac{\partial \det(M)}{\partial M_{ij}} = M^c_{ij}
\ea
where $M^c$ is the cofactor matrix of $M$.

Given a square matrix $M$ of nullity $1$, with a non-zero left kernel 
row vector $l$ and a non-zero right kernel vector $r$, its cofactor matrix
will be a rank-$1$ matrix with 
\bna
\label{eq:rl}
M^c \propto l^tr^t
\ena
(This follows from the fact that $M^tM^c= M^cM^t = \det(M)I$, and the 
fact that under our assumptions, there is at least one
non-zero cofactor).

To complete the chain rule, we see how each entry of $R(\hp)$ changes
with changes in $\hp$ giving us 
\ba
[d[\det(R(\hat{\p}))]]_{ik} 
&=& \sum_{j\in N(i)}  (R(\hat{\p})^c)^{\{ij\}}_{ik}
-
\sum_{j\in N_k(i)}
(R(\hat{\p})^c)^{\{ij\}}_{jk}
\ea
Like above,  $ik$ corresponds to an unpinned degree of freedom. 

Using the fact that $\omega$ is the left kernel of $R(\hat{\p})$
and $\p'$ is its right kernel, together with Equation (\ref{eq:rl}),
this gives us
\ba
[d[\det (R(\hat{\p}))]]_{ik} 
&=& \alpha\left[ \sum_{j\in N(i)}  \omega_{\{ij\}}(p'_i)_k - \sum_{j\in N_k(i)} \omega_{\{ij\}}(p'_j)_k \right]
\ea
where $\alpha$ is a global scale, 
proving the equivalence.

\paragraph{Generalizations}
In the previous section, we looked at the case of 
an isostatic graph with $m=nd-D$ edges, that was generically
rigid but had one infinitesimal flex at $\hp$. What about
other scenarios?

Let us first consider a generically rigid graph with more edges, say
$m=nd-D+1$, (hyperstatic), 
and a framework $\hp$  that has one (pinned) infinitesimal
flex $\p'$. In this case the pinned rigidity matrix has one more row
than column, and $R(\hp)$ will have a two dimensional space $S$
of stresses. 
In this case, the prestress stability certificate
is the existence of an $\omega \in S$ such that Equation~\ref{eq:fred} is satisfied. 
To search for such a certificate,
it is sufficient to check a basis of $S$.

Likewise, for a transversality test, we would look at all 
of the square matrices obtained by dropping one row from the
pinned $R(\hp)$ and certify that there is at least one such
square case where we find tranversality, a-la Equation (\ref{eq:trans}).
Similar to last section, this can be seen as  equivalent to finding a 
prestress stability certificate that based on a stress that has a 
zero entry on the edge corresponding to the dropped row.
Note, that by dropping each row from the rigidity matrix, we are also certain
to find a basis for the stress space $S$.

Next let us go back to an isostatic graph,
where $m=nd-D$, but now
suppose that the rank of $\hR(\hp) = nd-D-2$, ie it has
two  non-trivial flexes, and two  stresses.
In this case, there is still
a well defined notion of a prestress stability test. Indeed in
this case there is also a distinct notion of a second-order rigidity
test. 
Notably, 
in this case, $(R(\p))^c$ will be the all zero matrix, thus
$d[\det(R(\hp))]=0$ and so the transverse rigidity test 
of Equation (\ref{eq:trans})
will always fail!  

Finally, prestress stability and second order rigidity can be used to certify the rigidity 
of some framework $(G,\hp)$ where $G$ is generically flexible (say hypostatic).
It does not appear the transverse rigidity has any role in this setting, as $(G,\hp)$ must always 
be infinitesimally flexible. 

In summary, when $G$ is generically rigid, and $(G,\hp)$ has one infinitesimal flex, then 
transverse rigidity is equivalent to prestress stability. In all other cases, transverse
rigidity never holds.

\def\v{\oldv}
\bibliographystyle{abbrvlst}
\bibliography{framework}

\end{document}